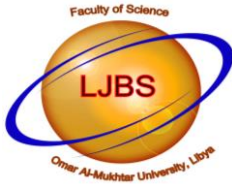
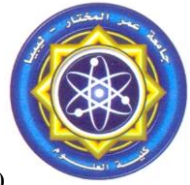

# Maximal Order of an NG-group

**Faraj. A. Abdunabi**

Department mathematical, Faculty science, University, Ajdabyia, Libya
Correspondence author: Faraj.a.abdunabi@uoa.edu.ly

## Abstract

This study was aimed to consider the NG-group that consisting of transformations on a nonempty set *A* has no bijection as its element. In addition, it tried to find the maximal order of these groups. It found the order of NG-group not greater than n. Our results proved by showing that any kind of NG-group in the theorem be isomorphic to a permutation group on a quotient set of *A* with respect to an equivalence relation on *A*.

**Keywords:** NG-group, permutation group, equivalence relation, $\chi$-subgroup.

## Introduction

This study considered the problem that the maximal order of a group consisting of transformations on a nonempty set *A* and the group has no bijection as its element.  Recall a permutation group on *A* is a group consisting of bijections from *A* to *A* with respect to compositions of mappings. It is well known that any permutation group on a set *A* with cardinality n has an order not greater than *n*!.

In previous studies, there are some authors (*1, 2*), problem 1.4 in (*3*) considering groups which consist of non-bijective transformations on *A* where the binary operation is the composition of mappings. Our first result is on the orders of such groups.

**Theorem1.1**. let *A* be a set with cardinality *n*. Suppose *NG be* groups consisting of non-bijective transformations on *A*, where the binary operation on *NG* is the composition of transformation. Then the order of *NG* is not greater than $(n-1)!$ and there are such groups having order $(n-1)!$.

Then it was proven Theorem1.1 by showing that any kind of group in the theorem be isomorphic to a permutation group on a quotient set of *A* with respect to an equivalence relation on *A*.

**Definition 1.1**. A class of group $\chi$ is called an SHP-class if it is closed under taking subgroups, homomorphic images, and products of normal subgroups. The latter condition means that if *U* and *V* are normal in *G* and both *U* and *V* lie in $\chi$, then $UV \in \chi$. If a group *G* belongs to $\chi$, we will say *G* is an $\chi$-group.

**Remark 1.2**. If $\chi$ is an SHP-class and $U,V \triangleleft G$ are such that $G/U$ and $G/V$ are $\chi$- groups, then $G/(U \cap V)$ is isomorphic to a subgroup of the $\chi$-group $(U/G) \times \{G/V\}$, and thus $G/(U \cap V)$ is an $\chi$-group. It follows that given a finite group *G*, there exists a unique smallest normal subgroup *N* such that $G/N \in \chi$, and we write $N = G\chi$.

The following theorem was found by the author; see also lemma 2.32 in (*4*).

**Theorem 1.2**. Let $\chi$ be an SHP-class, and suppose $G=UV$, where *U* and *V* are subnormal in *G*. then $G\chi = U\chi\, V\chi$.

It could take the SHP-class to the class of p-groups, the class of nilpotent groups, etc. Theorem1.2 will imply Lemma 9.15, problem 9B.5, Corollary 9.27, problem 9C.2, as corollaries.

**Remark 1.3.** It was noted in sec.4 of (*5*) that if it replaces the condition that $\chi$ is an SHP-class by some weaker condition that the class $\chi$ is such that whose composition factors all lie in some given set of simple groups then theorem 1.2 will fail in this case.

**Definition 1.2**. Let $\chi$ be an SHP-class and *G* be a finite group. The result was denoted the maximal normal $\chi$-subgroup of *G* by $G\chi$.

Then was considered the question that if $G=UV$ with *U,V* subnormal in *G* then it holds that $G\chi = U\chi V\chi$ or not. If *p* is a prime and take the SHP-class $\chi$ to be the class of all finite *p*-group, then for any finite group *G*χ will be $O_p(G)$ and results have the following theorem.

**Theorem1.3.** Let *p* and *q* be two primes such that $q \equiv 1 \pmod{p}$. Let $N=G_q$ be a cyclic group of order *q* and $H=<x>\times<y>$ an elementary abelian group of order $p^2$. Let $<x>$ act on *N* faithfully and $<y>$ act on *N* trivially. Set $G = N \rtimes H$ to the semidirect product of *N* and *H*. Let $U=N<x>$ and $V=<xy>$. Then

[1]  *U,V* are both subnormal in *G*.

[2] $O_p(G)=<y>$ and $O_p(U)=O_p(V)=1$. In particular, $Op(G) \neq O_p(U)O_p(V)$.

**Definition 2.1**. A binary relation ~ in *A* is called an equivalence relation on *A*. If it satisfies the following three conditions:
   (i)  *a~a* for any $a \in A$;
   (ii) for any $a,b \in A$, if *a~b* then *b~a* ;
   (iii) for any $a,b,c \in A$, if *a~b* and *b~c* then *a~c*.

For the set A, where use $A^A$ to dented the set of all its transforms, for any $f \in A^A$, we use Im(*f*) to denote the image of *f*. Also, *Z* and $Z_>$ will respective dented the set of integers and positive integers.

**Definition 2.2**. Let ~ be an equivalence relation on *A*, for an element $a \in A$, it is call $\{x \in A / x \sim \}$ the equivalence class of *a* determined by ~, which is denoted by $[a]_\sim$. And $A/\sim = \{[a]_\sim | a \in A\}$ is called the quotient set of *S* relative to the equivalence relation ~.

**Lemma2.3**. Theorem 1 (*1*),  for any $f \in G$ and the *e* the identity element of *G*, $\sim_e = \sim_f$

*Proof,* for any $a \in A$, so result goal is to show that
$[a]_f = [a]_e$.

On one hand, if $x \in [a]_f$, i.e. $f(x)=f(a)$. Since $G$ is a group with identity element $e$, there is a transformation $f' \in G$ such that $f'f=e=ff'$. Therefore,
$$e(x)=f'(f(x))=f'(f(a))=e(a),$$
Which yields that $x \in [a]_e$.

On the other hand, if $y \in [a]_e$ i.e. $e(a)e(y)$. Hence,
$$f(a)=(fe)(a)=f(e(y))=(fe)(y)=f(y),$$
Which implies $y \in \in [a]_f$. It follows that $[a]_e = [a]_f$ for any $a \in A$, as wanted.

**Remark 2.1.** For Lemma 2.4, the current result see that $\sim_f = \sim_g$ for any element $f, g \in G$

The following Theorem is the revised version of Theorem 2, (1).

**Theorem 2.5.** Let $f$ be an element in $A^A$ and $\hat{f}$ be the induced transformation of $f$ on $A/\sim_f$, i.e
$$\hat{f}: A/\sim_f \longrightarrow A/\sim_f,\ [x]_f \mapsto [f(x)]_f.$$
Then the following hold:

(i)   The exists a groups $G \subseteq A^A$ containing $f$ as the identity element iff $f^2 = f$.
(ii)  There is a groups $G \subseteq A^A$ containing $f$ as the identity element iff $\hat{f}$ is abijective on $A/\sim_f$.

The following two corollaries are from (1), and we make some corrections to the original proofs. Actually, this adopt the restriction of finiteness on $A$ in the first corollary from the original one. And then used the finiteness on $A$ in the second corollary; the original one did not use it.

**Corollary 2.6.** Let $f$ be an element in $A^A$. Then $f^2 = f$ iff the induced mapping $\hat{f}$ on $A/\sim_f$ is the identity element.

*Proof.* On one hand, suppose that $f^2 = f$. Then for any $[x]_f \in A/\sim_f$, as $f(x)=f(f(x))$, then see that $[x]_f = [f(x)]_f$. It follows that
$$\hat{f}([x]_f) = [f(x)]_f = [x]_f\ ;$$
This implies that $\hat{f}$ is the identity mapping on $A/\sim_f$.

On the other hand, assume that $\hat{f}$ is the identity mapping on $A/\sim_f$. Then for any $[x]_f \in A/\sim_f$, the condition that $\hat{f}([x]_f) = [x]_f$ will imply that $[f(x)]_f = [x]_f$ and hence $f(f(x)) = f(x)$. It follows that $f^2 = f$ as required.

**Corollary 2.7.** Suppose that $A$ is a finite set and $f$ is an element in $A^A$. Then there is a group $G \subseteq A^A$ containing $f$ as an element iff $\mathrm{Im}(f) = \mathrm{Im}(f^2)$.

*Proof*. On one hand, suppose that there is a group $G \subseteq A^A$ containing $f$ as an element. Let $e$ be the identity element of $G$. Then by Theorem 2.5, the induced mapping $\hat{f}$ is a bijection on $A/\sim_f$ . In particular, $\hat{f}$ is surjective and thus for any $x \in A$, there is a $[y]_f \in A/\sim_f$ such that $\hat{f}([y]_f) = [x]_f = [f(y)]_f$ ; Which yields that $f(x) = f(f(y)) = (f^2)(y)$. As a result, $\text{Im}(f) \subseteq \text{Im}(f^2)$ and thus $\text{Im}(f) = \text{Im}(f^2)$.

On the other hand, suppose that $\text{Im}(f) = \text{Im}(f^2)$. Thus, for any $f(x) \in \text{Im}(f)$ there is a $y \in A$ such that $f(x) = f(f(y))$ and hence $\hat{f}([y]f) = [x]f$ ; which implies that $\hat{f}$ is surjective on $A/\sim_f$ . Note that results are assuming that A is finite and so is $A/\sim_f$ . this study has that the induced mapping $\hat{f}$ is bijective. By Theorem 2. 5, the assertion follows.

**Remark** 2. 2. Let $G \subseteq A^A$ be a group. That has seen, in Remark 2.1, that $\sim_f = \sim_g$ for any elements in $G$ and we will denote the common equivalence relation by $\sim$. Also, by Theorem 2.5, each element $f \in G$ will induce a bijection $\hat{f}$ on $A/\sim$.

The following theorem is crucial since it turns a group $G \subseteq A^A$ into a permutation group.

**Theorem 2.8**. Let $G \subseteq A^A$ be a group. Set
$\hat{G} = \{\hat{f} \mid f \in G\}$; then $\hat{G}$ is a permutation group on $A/\sim$ and $\rho : G \to \hat{G}$, $f \mapsto \hat{f}$, is an isomorphism.

*Proof*. For any $f, g \in G$ and any $[a] \in A/\sim$, results have $\rho(fg)([a]) = [(fg)(a)] = [f(g(a))] = \rho(f)([g(a)]) = (\rho(f) \rho(g))([a])$; which implies that $\rho(fg) = \rho(f) \rho(g)$ and thus $\rho$ is a homomorphism. By the definition of $\hat{G}$, it is obvious that $\rho$ is surjective.
Now suppose that $\rho(f) = \rho(g)$ for two elements $f, g \in G$, i.e.$[f(a)] = [g(a)]$, $\forall a \in A$: Let $e$ be the identity element of $G$, then we have $[f(a)]_e = [g(a)]_e$; $\forall a \in A$. It follows that $e(f(a)) = e(g(a))$; $\forall a \in A$ .Hence,$f(a) = (ef)(a) = e(f(a)) = e(g(a)) = g(a)$, $\forall a \in A$, and therefore $f = g$. so it conclude that $\rho$ is injective. As a consequence, $\rho$ is an isomorphism.

**Definition 2.3**. A subgroup $H$ of a group $G$ is called characteristic in $G$, denoted $H$ char $G$, if every automorphism of $G$ maps $H$ to itself, that is $\rho(H) = H$ for all $\rho \in \text{Aut}(G)$.

**Remark 2.3**. If $H$ is characteristic in $G$ in $K$ and $K$ is characteristic in $G$, then $H$ is characteristic in $G$.

Let $G$ be a finite group. It has the following two lemmas. They are from Section 2 of (5).

**Lemma 2.9**. Suppose that $\chi$ is an SHP-class.

(a) Let $\leq G$ be a subgroup. Then $H^\chi \leq G^\chi$.
(b) Let $N \triangleleft G$ be a normal subgroup of $G$ and write $\overline{G}=G/N$. then $\overline{G}^\chi = \overline{G^\chi}$.
(c) $G^\chi$ is characteristic in $G$.
(d) $O_\chi(G)$ is characteristic in $G$.

The following lemma is a generalization of Problem 2A.1 in (8).

**Lemma 2.10**. Let $A$ and $B$ be two subnormal $\chi$-subgroups of $G$. Then the subgroup $<A,B>$ generated by $A$ and $B$ are $\chi$-subgroup of $G$.

*Proof*. Let $A$ be a subnormal $\chi$-subgroup of $G$. The resulting use induction on the subnormal depth $r$, $A \subseteq O_\chi$ of $A$ in $G$ to show that if $r = 1$, then $A \triangleleft$ and thus $A \subseteq O_\chi(G)$ since $O_\chi(G)$ is the largest normal $\chi$-subgroup of $G$.

Suppose $r > 1$ and the containment holds for r-1. Let $A_1 = A \triangleleft ... \triangleleft H_{r-1} \triangleleft H_r = G$ be a subnormal series from $A$ to $G$: Then $A \subseteq O_\chi(G)$ by inductive hypothesis. Since $O_\chi(G)$ char $H_{r-1}$ and $H_{r-1} \triangleleft G$; $O_\chi(G) \triangleleft G$ and then $O_\chi(G) (H_{r-1}) \subseteq O_\chi(G)$. It was concluded that $A \subseteq O_\chi(G)$.
In general, for any two subnormal $\chi$-subgroups $A$ and $B$, $A, B \subseteq O_\chi(G)$ and thus $<A,B> \subseteq O_\chi(G)$ as wanted.

## Proofs of Main Results

Now let $A$ be a set having $n$ letters written as $\{1, 2, ..., n\}$. The results have the following theorem, which is Theorem 1.1.

**Theorem 3.1**. Let $A$ be a set with cardinality n with $n \geq 3$. Suppose $NG$ is a group consisting of non-bijective transformations on $A$, where the binary operation on $NG$ is the composition of transformations. Then the order of $NG$ is not greater than $(n-1)!$ and there are such groups having order $(n-1)!$:

*Proof*. Let $NG$ be a group consisting of non-bijective transformations on $A$. By Remark 2.1, it is known that $\sim_f = \sim_g$ for any element $f, g \in NG$ and it denote the common equivalence relation by $\sim$. Note that $NG$ is a group consisting of non-bijective transformations, then we see that the equivalence relation is not the equality relation = on $A$. Thus, these results have that the quotient set $A/\sim$ has an order less than $n-1$.

Additionally, *NG* is isomorphic to a permutation group on *A/~* by Theorem 2.8. It follows that the order of *NG* is less than (*n*-1)! as any permutation group on *A/~* has order less than (*n*-1)!.

Note that In defining a permutation s on the set {1, 3, …, *n*}, there are *n*-1 choices for $\rho(1)$, *n*-2 choices of $\rho(3) \neq \rho(1)$, *n*- 2 choices of $\rho(4) (\neq \rho(1), \rho(3))$, etc., i.e. totally (*n*-1)(*n*-2)1 = (*n*-1)!.

**Theorem 3.2**. Let $\chi$ be an SHP-class, and suppose $G = UV$; where *U* and *V* are subnormal in *G*. Then $G^\chi = U^\chi V^\chi$.

*Proof*. This work use induction of the subnormal depth of *U* in *G* to prove the result.

First, if the subnormal depth of *U* in *G* is one, i.e. $U \triangleleft G$. Since $U^\chi$ is characteristic in *U* and *U* is normal in *G* we see that $U^\chi$ is normal in *G*. Let $\bar{G} = G/U^\chi$. By the hypothesis, $\bar{G} = \bar{U}\bar{V}$ where $\bar{U} = U/U^\chi$, $\bar{V} = VU^\chi/U^\chi$. Thus, $\bar{U}$ is a normal $\chi$-group of $\bar{G}$ and $\bar{V}$ is subnormal in $\bar{G}$. By Lemma 2.10, we have $\bar{G}^\chi = \bar{V}^\chi$. By Lemma 2.9 (b), $\bar{G}^\chi = \overline{G^\chi}$, $\bar{V}^\chi = \overline{V^\chi} = \overline{U^\chi V^\chi}$.
By correspondence theorem, it has $G^\chi = U^\chi V^\chi$; as required.

Now suppose that the subnormal depth of *U* in *G* is r with $r > 1$: Let
$U_1 = U \triangleleft \ldots \triangleleft U_r \triangleleft G$
be a subnormal series from *U* to *G* with length *r*. By Dedekind's lemma,
$U_r = U(V \cap U_r)$. As both *U* and $V \cap U_r$ are subnormal in $U_r$ and *U* has subnormal depth *r*- 1 in $U_r$, then obtain that
$(U_r)^\chi = U^\chi (V \cap U_r)^\chi$
by inductive hypothesis. Also, $G = U_r V$ with $U_r$ normal in *G* and *V* subnormal in *G*, and hence
$G^\chi = (U_r)^\chi V^\chi$ by the first paragraph of the proof. It follows that
$G^\chi = (U_r)^\chi V^\chi = U^\chi (V \cap U_r)^\chi V^\chi = U^\chi V^\chi$,
because $(V \cap U_r)^\chi \subseteq V^\chi$ by Lemma 2.9 (a).

**Theorem 3.3**. Let *p* and *q* be two primes such that $q \equiv 1 \pmod{p}$. Let $N = C_q$ be a cyclic group of order q and $H = \langle x \rangle$ an elementary abelian group of order $p^2$. Let $\langle x \rangle$ act on *N* faithfully and $\langle y \rangle$ act on *N* trivially. Set $G = N \rtimes H$ to be the semidirect product of *N* and *H*. Let $U = N \langle x \rangle$ and $V = N \langle xy \rangle$. Then
  (i)  *U*, *V* are both subnormal in *G* and $G = UV$.
  (ii) $O_p(G) = \langle y \rangle$ and $O_p(U) = O_p(V) = 1$. In particular, $Op(G) \neq O_p(U)O_p(V)$.

*Proof*. Since *N* is normal in *G* and the quotient group *G=N/H* is abelian, it deduced that the derived subgroup *G'* is contained in *N*. It follows that both *U* and *V* contain *G'* as a subgroup, which implies that *U* and *V* are normal in *G*. Obviously, *G = UV*. Assertion (i) holds.

Note that the Sylow *p*-subgroup of *G* is not normal since $<x>$ act on *N* faithfully and hence $O_p(G)$ has an order less than $p^2$. However, as $<y>$ act on *N* trivially, *N* normalizes $<y>$ which yields that $<y>$ is a normal *p*-subgroup of *G*. It is easy to see that $O_p(G) = <y>$. Both $<x>$ and $<xy>$ act faithfully on *N*, which yields that $O_p(U) = O_p(V) = 1$; as wanted.

# ترتيب لزمر NG تبية

## فرج أرخيص عبدالنبي

قسم الرياضيات, كلية العلوم, جامعة اجدابيا

## الملخص العربي

هدفت هذه الدراسة إلى اعتبار زمر NG التي تتكون من تحولات غير تقابليه على مجموعة غير فارغة A والتي تكون زمر غير جزئية من الزمر التبديلية ونطلق على هذه الزمر بـ NG-Transformation. وسنتحصل علي انه لايمكن ان يكون درجة اعلي زمرة أكبر من (ن-1)!. بالإضافة إلى ذلك، سوف نثبت نتيجتنا من خلال إظهار أن أي نوع من هذه الزمر في النظرية المتحصل عليها سيكون متماثلًا لزمرة التباديل المعرفة علي مجموعة القسمة A بالنسبة لعلاقة التكافؤ على A.

الكلمات المفتاحية: مجموعة NG ، زمر تبديلية ، علاقة تكافؤ، المجموعة الفرعية x.